\documentclass[a4paper]{article}

\usepackage{amssymb} 
\usepackage{amsmath} 
\usepackage{latexsym} 
\usepackage{theorem} 

\usepackage{makeidx} 

\usepackage{epsfig}

\def\a{\alpha}
\def\de{\delta}
\def\e{\varepsilon}
\def\la{\lambda}
\def\La{\varLambda}
\def\p{\prime}
\def\w{\wedge}
\def\vect#1{\mbox{\boldmath $#1$}} 
\def\smallvect#1{\mbox{\small{\boldmath $#1$}}}
\def\Eap{E^{\,\a,p}}
\def\omegasub#1{\mbox{\large $\omega$}\mbox{\small${\mbox{\large ${}$}}_{#1}$}}%
\theorembodyfont{\itshape} 
\theoremstyle{break} 

\def\RR{\mathbb{R}}

\def\NN{\mathbb N}

\newtheorem{theorem}{Theorem}
\newtheorem{lemma}[theorem]{Lemma}

\newtheorem{proposition}[theorem]{Proposition}

\theorembodyfont{\rmfamily} 
\theoremstyle{change} 

\numberwithin{equation}{section}

\begin{document}

\title{\bf Errata and Comments for \\ ``Energy of knots and conformal geometry''}

\author{Jun ~O'Hara 
\\{}
\\{\small Department of Mathematics and Informatics,Faculty of Science, 
Chiba University}
\\{\small e-mail : ohara@math.s.chiba-u.ac.jp
}}

\maketitle


\begin{abstract}
This article serves as errata of the book 

\smallskip
\begin{center}
``{\em Energy of knots and conformal geometry}", Series on Knots and Everthing Vol. 33, 

World Scientific, Singapore, 304 pages, (2003). 
\end{center}

\smallskip
The updated version is available through web, linked from the author's Home Page:

\smallskip
\centerline{
http://www.comp.tmu.ac.jp/knotNRG/indices/indexe.html
}

\medskip
The list below would be far from being completed\footnote{The author can find errors almost every time he has a look at the book.}. 

Suggestions and comments would be deeply appreciated. 
\end{abstract}

\bigskip
{\small {\it Key words and phrases}. Energy, knot, link, M\"obius geometry, conformal geometry, cross ratio}

{\small 1991 {\it Mathematics Subject Classification.} Primary 57M25, Secondary 53A30}

\begin{itemize}

\item Throughout the book ${}_{n}C_{p}$ means the number of combination ${n\choose{p}}=\frac{n!}{p!(n-p)!}$. The author is sorry for the Japanese notation. 

\item Page 3. The second paragraph 

``Let $M$ denote $\mathbb{R}^3$ or $S^3$. 
Two knots $f$ and $f^{\prime}$ in $M$ are called {\bf isotopic}
 if there is an isotopy $h_t:M\to M$ $(t\in [0,1])$ of the ambient space
 such that $h_0$ is equal to the identity map and that the map
 $(x,t)\mapsto (h_t(x),t)$ from $M\times [0,1]$ to itself is a homeomorphism. 
Then two knots are isotopic if and only if there is an orientation preserving
 homeomorphism $h$ of $M$ that satisfies $f^{\prime}=h\circ f$. 
A {\bf knot type} $[K]$ a knot $K$ is an isotopy class of $K$. "

should be replaced by 

``Let $M$ denote $\mathbb{R}^3$ or $S^3$. 
Two knots $f$ and $f^{\prime}$ in $M$ are called {\bf isotopic}
 if there is an isotopy $h_t:M\to M$ $(t\in [0,1])$ of the ambient space
 such that $h_0$ is equal to the identity map, that the map
 $(x,t)\mapsto (h_t(x),t)$ from $M\times [0,1]$ to itself is a homeomorphism, {\bf and that \boldmath $h_1\circ f=f^{\prime}$}. 
Then two knots are isotopic if and only if there is an orientation preserving
 homeomorphism $h$ of $M$ that satisfies $f^{\prime}=h\circ f$. 
{\bf In this book, let us call an isotopy class of a knot $K$ a {\bf knot type} of $[K]$. }"

\item Page 15. The second line from the bottom. 
$$E_{\e}^{(\a)}(K)=E_{\e}^{(\a)}(h)=\int_{\vect{S^1}}V_{\e}^{(\a)}(K;x)dx
=\int_{S^1}V_{\e}^{(\a)}(h;s)ds.$$
should be replaced by 
$$E_{\e}^{(\a)}(K)=E_{\e}^{(\a)}(h)=\int_{\vect{K}}V_{\e}^{(\a)}(K;x)dx
=\int_{S^1}V_{\e}^{(\a)}(h;s)ds.$$

\item Page 26. ``then $E^{(\a )}$" in the first line of Theorem 2.4.1 (2) should be removed. 

\item Page 26. The 5th line from the bottom, i.e. the second line in the proof of Theorem 2.4.1 (2). 

``Let $h$ be a knot with $\vect{|h^{\prime\prime}|\equiv 1}$ and $b=E^{(\a)}(h)$ $(b>0)$. "

should be replaced by 

``Let $h$ be a knot with $\vect{|h^{\prime}|\equiv 1}$ and $b=E^{(\a)}(h)$ $(b>0)$. "

\item Page 34. {\bf Definition 3.1}. 
The first line of the formula holds if $X\ne P$: 
$$
I_{\varSigma}(X)=\left\{\begin{array}{ll}
\displaystyle{P+\frac{r^2}{|X-P|^2}(X-P)} & \quad \mbox{if} \quad X\in\mathbb{R}^3,\\
\infty & \quad \mbox{if} \quad X=P,\\[2mm]
P & \quad \mbox{if} \quad X=\infty.
\end{array}\right. 
$$
should be replaced by 
$$
I_{\varSigma}(X)=\left\{\begin{array}{ll}
\displaystyle{P+\frac{r^2}{|X-P|^2}(X-P)} & \quad \mbox{if} \quad X\in\mathbb{R}^3{\vect{\setminus\{P\}}},\\
\infty & \quad \mbox{if} \quad X=P,\\[2mm]
P & \quad \mbox{if} \quad X=\infty.
\end{array}\right. 
$$

\item Page 42. Outline of proof of Theorem 3.5.1. 

The convergence means the convergence with respect to the $C^0$-topology. 

\item Page 43. The first line. 
``Kusner and J. Sullivan'' should be replaced by ``Kim and Kusner''. 

\item The book does not contain the proofs of Theorem 3.6.1 on page 44 and 
Theorem 3.7.7 on page 53, 
for which the reader is referred to [FHW] and [He1] respectively. 

\item Page 45. The right hand side of the 6th line from the bottom. 

$$-\e\left[(1-\xi)\left\{\vect{{|f^{\prime}|}^{\prime}}(s+\e\xi)\right\}\right]_0^1+\e^2\int_0^1(1-\xi)\left\{{|f^{\prime}|}^{\prime}(s+\e\xi)\right\}\,d\xi$$

should be replaced by 

$$-\e\left[(1-\xi)\left\{\vect{{|f^{\prime}|}}(s+\e\xi)\right\}\right]_0^1+\e^2\int_0^1(1-\xi)\left\{{|f^{\prime}|}^{\prime}(s+\e\xi)\right\}\,d\xi$$

\item Page 54. Section 3.8. 

The main point of [KK] is as follows: 
As $r$ approaches $0$ or $1$, the orbital configurations approach a $p$- or $q$-fold covered circle, and thus its energy $E=E^{(2)}$ must go to infinity. 
Since $E$ is continuous, it takes on a minimum at some $r_0$, and by the symmetric criticality argument, this is actually critical among all variations, not just the orbital ones. 

\item Page 57. The 4th line from the bottom. 
$$
\int_{|z|=1}\left(
\frac{i(r^2p^2+(1-r^2)q^2)}{r^2z^{1-p}(z^{p}-1)^2
+(1-r^2)z^{1-q}(z^{q}-1)^2}
-\frac{i}{(z-1)^2}\right)\vect{ds}
$$
should be replaced by
$$
\int_{|z|=1}\left(
\frac{i(r^2p^2+(1-r^2)q^2)}{r^2z^{1-p}(z^{p}-1)^2
+(1-r^2)z^{1-q}(z^{q}-1)^2}
-\frac{i}{(z-1)^2}\right)\vect{dz}
$$

\item Page 58, The 8th line in subsection 3.9.1. 
$$\begin{array}{c}
{\rm Conf}_{n,m}(K,\vect{S^3})=\left.\left\{(x_1, \cdots, x_{n+m})\in K^n\times (\vect{S^3})^m\,\right|\,
x_j\ne x_k \, (j\ne k)\right\}\,,\\[2mm]
U_D=\left.\left\{(x_1, \cdots, x_{n+m})\in{\rm Conf}_{n,m}(K,\vect{S^3})\,\right|\,
x_1\prec \cdots \prec x_n\right\}\,.
\end{array}
$$
should be replaced by 
$$\begin{array}{c}
{\rm Conf}_{n,m}(K,\mathbb{R}^3)=\left.\left\{(x_1, \cdots, x_{n+m})\in K^n\times (\mathbb{R}^3)^m\,\right|\,
x_j\ne x_k \, (j\ne k)\right\}\,,\\[2mm]
U_D=\left.\left\{(x_1, \cdots, x_{n+m})\in{\rm Conf}_{n,m}(K,\mathbb{R}^3)\,\right|\,
x_1\prec \cdots \prec x_n\right\}\,.
\end{array},
$$
i.e. $S^3$ should be replaced by $\mathbb{R}^3$, as we use Euclidean metric in this subsection. 

\item Page 59, in Definition 3.5. 
$$\begin{array}{l}
E_{X,\cos}(K)=\displaystyle{\int_{K^4;\,x_1\prec \smallvect{x_3} \prec \smallvect{x_2}\prec x_4}\,
\frac{(1-\cos\theta_{13})(1-\cos\theta_{24})}{|x_1-x_3|^2|x_2-x_4|^2}
\,dx_1dx_2dx_3dx_4,}\\[4mm]
E_{X,\sin}(K)=\displaystyle{\int_{K^4;\,x_1\prec \smallvect{x_3} \prec \smallvect{x_2}\prec x_4}\,
\frac{\sin\theta_{13}\sin\theta_{24}}{|x_1-x_3|^2|x_2-x_4|^2}
\,dx_1dx_2dx_3dx_4,}\end{array}$$
should be replaced by 
$$\begin{array}{l}
E_{X,\cos}(K)=\displaystyle{\int_{K^4;\,x_1\prec \smallvect{x_2} \prec \smallvect{x_3}\prec x_4}\,
\frac{(1-\cos\theta_{13})(1-\cos\theta_{24})}{|x_1-x_3|^2|x_2-x_4|^2}
\,dx_1dx_2dx_3dx_4,}\\[4mm]
E_{X,\sin}(K)=\displaystyle{\int_{K^4;\,x_1\prec \smallvect{x_2} \prec \smallvect{x_3}\prec x_4}\,
\frac{\sin\theta_{13}\sin\theta_{24}}{|x_1-x_3|^2|x_2-x_4|^2}
\,dx_1dx_2dx_3dx_4,}\end{array},$$
i.e. the order of $x_2$ and $x_3$ in ``$x_1\prec x_3 \prec x_2\prec x_4$" 
should be reversed. 

\item Page 60. Section 3.11. 
Let me give an idea of the proof of self-repulsiveness of surface energies. 

(1) {\bf Kusner-Sullivan's $\vect{(1-\cos\theta)^2}$ energy $\vect{E_{KS}}$}

Suppose a $2$-dimensional surface $M$ in $\RR^3$ (the dimension of the ambient space does not matter) is tangent to itself at a point $p_0$. 
Assume that the intersection of $M$ and a small neighbourhood of $p_0$ consists of two connected components, say $S_1$ and $S_2$, each of which is almost flat. 
Assume that $S_1$ and $S_2$ can be expressed as graphs of functions $f_1$ and $f_2$ defined on a subset $W$ of the common tangent space $\Pi:=T_{p_0}S_1=T_{p_0}S_2$. We use coordinates of $\Pi$ so that $p_0$ is the origin. 

\begin{proposition}
We can find a positive constant $C$ and a small region $U$ in $W$ so that 
\begin{itemize}
\item Put $U_i:=2^{-i}U$ $(i\in\overline{\NN}:=\NN\cup\{0\})$, then $U_i\cap U_j=\emptyset$ if $i\ne j$. 
\item Put $S^1_i:=f_1(U_i)$ and $S^2_i:=f_2(U_i)$, then 
\[
\int_{S^1_i}\int_{S^2_i}\frac{(1-\cos\theta)^2}{|x-y|^4}\,dxdy\ge C \hspace{0.4cm}(\forall i\in\overline{\NN}).
\]
\end{itemize}
\end{proposition}

The claim implies that 
\[
E_{KS}(M)\ge\sum_{i=0}^{\infty}\int_{S^1_i}\int_{S^2_i}\frac{(1-\cos\theta)^2}{|x-y|^4}\,dxdy=\infty,
\]
which implies the self-repulsiveness of $E_{KS}$. 
Here, we used the continuity of $E_{KS}$ with respect to $C^2$-topology, which could be proved anyhow. 

The claim would be proved in the following way. 

(1) The dominat term of the integral can be estimated by up to the quadratic term of $f_1,f_2$. 

(2) If $x_i\in S^1_i$ and $y_i\in S^2_i$ then 
\begin{equation}\label{f1}
|x_i|\sim|y_i|\sim 2^{-i}, \hspace{0.2cm} |x_i-y_i|\sim |x_i|^2,  \hspace{0.2cm} |\theta(x_i,y_i)|\sim|x_i|,  \hspace{0.2cm} 1-\cos\theta(x_i,y_i)\sim|x_i|^2, 
\end{equation}
where $A_i\sim B_i$ means there are positive constants $C'$ and $C''$ that are independent of $i$ such that $C'A_i\le B_i \le C''A_i$ for all $i$. 
Since $\mbox{Area}(S^1_i)\sim\mbox{Area}(S^2_i)\sim|x_i|^2$ we have
\[
\int_{S^1_i}\int_{S^2_i}\frac{(1-\cos\theta)^2}{|x-y|^4}\,dxdy \sim 1, 
\]
which shows the proposition. 

\medskip
(2) {\bf Auckly-Sadun's regularized energy $\vect{E_{AS}}$}

Let $\mathcal{S}_C$ be a set of closed surfaces $M$ in $\RR^3$ of class $C^5$ that are ``{\sl bounded by a positive constant $C$ in the sense of $C^5$-topology}'', to be precise, 
for any point $x$ on $M$, if we express a small neighbourhood of $x$ of $M$ as a graph of a function $f$ on a small neighbourhood $W_x$ of the origin of the tangent plane $T_xM$, then $\left|\partial^\a f\right|\le C$ on $W_x$ for any multi-index $\a=(\a_1,\a_2)$ with $0\le |\a|=\a_1+\a_2\le 5$. 

Assume $M\in \mathcal{S}_C$. 
Let $k_{M}$ be the maximum of the principal curvatures of $M$. Note that $k_{M}\le C$. Put $R_0:=1/C$, then $R_0\ge1/ k_{M}$ for any $M\in \mathcal{S}_C$. 

For $x\in M$ and $0<r\le R_0$, let $U_r(x)$ be the connected component of $M\cap B_r(x)$ that contains $x$, where $B_r(x)$ is the $3$-ball with center $x$ and radius $r$. 

Recall that 
\[
E_{AS}(M):=\int_M \left[\lim_{\e\to0}\left(\int_{M, |x-y|\ge\e}\frac{dy}{|x-y|^4}-\frac\pi{\e^2}+\frac{\pi\Delta(x)}{16}\log\left(\Delta(x)\e^2\right)+\frac{\pi K(x)}4\right)\right] dx,
\]
where $\Delta(x):=(\kappa_1-\kappa_2)^2$ and $K(x)=\kappa_1\kappa_2$, where $\kappa_1$ and $\kappa_2$ are principal curvatures of $M$ at $x$. 
Put 
\[
\begin{array}{rcl}
V(U_r(x);x)&:=&\displaystyle \lim_{\e\to0}\left(\int_{U_r(x), |x-y|\ge\e}\frac{dy}{|x-y|^4}-\frac\pi{\e^2}+\frac{\pi\Delta(x)}{16}\log\left(\Delta(x)\e^2\right)+\frac{\pi K(x)}4\right), 
\end{array}
\]
then 
\[
E_{AS}(M)=\int_M V(U_r(x);x)\,dx+\int_M\left(\int_{M\setminus U_r(x)}\frac{dy}{|x-y|^4}\right)dx. 
\]

\begin{lemma}
\begin{enumerate} 
\item There is a constant $C'(r)$ (which might be negative) such that 
$V(U_r(x);x)\ge C'(r)$ for any $M\in \mathcal{S}_C$, for any $x\in M$ and for any $r$ with $0<r\le R_0$. 

(I think we can take $C'(r)$ to be independent of $r$.)
\item There is a positive constant $A$ such that $\mbox{\rm Length}(\partial U_r(x))\ge Ar$ for any $M\in \mathcal{S}_C$, for any $x\in M$ and for any $r$ with $0<r\le R_0$. 
\end{enumerate}
\end{lemma}

Suppose $M\in\mathcal{S}_C$ is close to have a double point. We look for a lower bound of $E_{AS}(M)$ if $M$ satisfies satisfies 
\[
M\in\mathcal{S}_C, \>\> \exists\, p,q\in M, \>\> 0<\exists \, r\le R_0, \>\> U_r(p)\cap U_r(q)=\emptyset, \>\> |p-q|=d, 
\]
for some positive number $d$, and consider the limit as $d\downarrow0$. 

First note 
\[
E_{AS}(M)\ge C'\left(\frac{r}2\right)\textrm{Area}(M)+\int_{U_{r/2}(p)}\int_{U_{r/2}(q)}\frac{dx\,dy}{|x-y|^4}. 
\]
The second term of the right hand side is bounded below by 
\[
\int_0^{r/2}\int_0^{r/2}\frac{AsAt}{(d+s+t)^4}\,dsdt, 
\]
which blows up to $+\infty$ as $d$ goes down to $0$, which proves the $C^5$-self-repulsiveness of Auckly-Sadun surface energy.

\item Page 67. The 5th line 

If $\Eap (h)$ with $\a p>2$ is {\bf finite} then $h$ cannot have a sharp turn. 

 should be replaced by 

If $\Eap (h)$ with $\a p>2$ is {\bf bounded} then $h$ cannot have a sharp turn. 
\item Pages 95, Figure 6.7, page 97, Figure 6.11, page 99, Figure 6.15. 
The caption 
``Look with the {\bf right} eye." 
should be 
``Look with the {\bf left} eye." 
(You do not have to close your left eye.) 

\item There is a misunderstanding about the order of contact. 
The order of contact in the book should be reduced by $1$. 
The errors can be found on pages 119, 120, 123-125, 184-185. 

\begin{itemize}
\item Page 119. Definition 8.2 

``(1) An {\bf osculating circle} $\cdots$ 
is the circle which is tangent to $K$ at $x$ at least to the {\bf third} order''

\smallskip
should be replaced by 
\smallskip

``(1) An {\bf osculating circle} $\cdots$ 
is the circle which is tangent to $K$ at $x$ at least to the {\bf second} order''

and 

``(2) An {\bf osculating sphere} $\cdots$ is tangent to $K$ at $x$ at least to the {\bf fourth} order. ''

should be replaced by 

``(2) An {\bf osculating sphere} $\cdots$ is tangent to $K$ at $x$ at least to the {\bf third} order. ''

\medskip
\item Page 119. 

``{\rm \bf Proposition 8.3.1}{\it 
{\rm (1)} An osculating sphere is uniquely determined 
if the order of tangency of the osculating circle to the knot is just $\vect 3$. 

\smallskip
{\rm (2)} Suppose the order of tangency of the osculating circle 
to the knot is just $\vect 3$. 
Then $\cdots$ }''

\smallskip
should be replaced by 
\smallskip

``{\rm \bf Proposition 8.3.1}{\it 
{\rm (1)} An osculating sphere is uniquely determined 
if the order of tangency of the osculating circle to the knot is just $\vect 2$. 

\smallskip
{\rm (2)} Suppose the order of tangency of the osculating circle 
to the knot is just $\vect 2$. 
Then $\cdots$ }''

\medskip
\item Page 120. After Proposition 8.3.1. 

``When the order of tangency of the osculating circle to the knot 
is greater than $\vect 3$, any $2$-sphere through the osculating circle 
is tangent to the knot to the {\bf fourth} order. 
But there might be a unique $2$-sphere 
which is tangent to the knot with a higher order of tangency than $\vect 4$. ''

\smallskip
should be replaced by 
\smallskip

``When the order of tangency of the osculating circle to the knot 
is greater than $\vect 2$, any $2$-sphere through the osculating circle 
is tangent to the knot to the {\bf third} order. 
But there might be a unique $2$-sphere 
which is tangent to the knot with a higher order of tangency than $\vect 3$. ''

\medskip
\item Page 120. The 6th line from the bottom (after the formula (8.1)). 

``Then $C$ is tangent to $K$ at $0$ to the {\bf fourth} order, 
i.e. $f^{(3)}(0)=C^{(3)}(0)$, if and only if $k^{\p}=0$ and $k\tau=0$. ''

\smallskip
should be replaced by 
\smallskip

``Then $C$ is tangent to $K$ at $0$ to the {\bf third} order, 
i.e. $f^{(3)}(0)=C^{(3)}(0)$, if and only if $k^{\p}=0$ and $k\tau=0$. ''

\medskip
\item Pages 123-124. Proof of Proposition 8.4.2. 

``Case I: $\cdots$

If the knot $K$ is transversal to $\varSigma$ at $x$ and if 
the order of tangency of $K$ and $\varSigma$ at $y$ is $\vect 2$, 
then $K$ must intersect $\varSigma$ (not necessarily transversally) 
at a third point $z$ which is different from both $x$ and $y$, 
and therefore $\varSigma=\sigma(x,z,y,y)$. 

If the order of tangency of $K$ and $\varSigma$ at $y$ is more than 
or equal to $\vect 3$, then $\varSigma=\sigma(x,y,y,y)$. 

\noindent Case II: Suppose $x=y$. 
The order of tangency of $K$ to $\varSigma$ at $x$ is more than or equal to $\vect 3$. 

If it is $\vect 3$ then $K$ must intersect $\varSigma$ (not necessarily transversally) 
at another point $z$ $(z\ne x)$. 
Then $\varSigma=\sigma(x,x,x,z)$. 

If the order of tangency is more than or equal to $\vect 4$ then 
$\varSigma=\sigma(x,x,x,x)$. ''

\smallskip
should be replaced by 
\smallskip

``Case I: $\cdots$

If the knot $K$ is transversal to $\varSigma$ at $x$ and if 
the order of contact of $K$ and $\varSigma$ at $y$ is $\vect 1$, 
then $K$ must intersect $\varSigma$ (not necessarily transversally) 
at a third point $z$ which is different from both $x$ and $y$, 
and therefore $\varSigma=\sigma(x,z,y,y)$. 

If the order of contact of $K$ and $\varSigma$ at $y$ is more than 
or equal to $\vect 2$, then $\varSigma=\sigma(x,y,y,y)$. 

\noindent Case II: Suppose $x=y$. 
The order of tangency of $K$ to $\varSigma$ at $x$ is more than or equal to $\vect 2$. 

If it is $\vect 3$ then $K$ must intersect $\varSigma$ (not necessarily transversally) 
at another point $z$ $(z\ne x)$. 
Then $\varSigma=\sigma(x,x,x,z)$. 

If the order of contact is more than or equal to $\vect 3$ then 
$\varSigma=\sigma(x,x,x,x)$. ''

\item Page 125. The 1st line. 

``Since $\rho(x)>r(\varSigma)$ by the assumption, the knot $K$ cannot have 
the tangency of order $\vect 3$ with $\varSigma$ at $x$, and hence $K$ must lie in $\cdots$''

\smallskip
should be replaced by 
\smallskip

``Since $\rho(x)>r(\varSigma)$ by the assumption, the knot $K$ cannot have 
the tangency of order $\vect 2$ with $\varSigma$ at $x$, and hence $K$ must lie in $\cdots$''

\end{itemize}

\item Page 125. The 5th line. 

``then 
the osculating sphere at $x$ contains the osculating {\bf sphere} at $x$ 
as the great circle''

should be replaced by 

``then 
the osculating sphere at $x$ contains the osculating {\bf circle} at $x$ 
as the great circle''

\item Page 141. The 7th and 8th lines from the bottom, i.e. the last two lines of (4) have two errors: 

Since $u_5$\underline{, }$v_5>1$ this implies 
$u_1v_1+ \cdots +u_4v_4\le u_5v_5-1$ which means $\langle \vect u,\vect v\rangle <-1$. 

should be replaced by 

Since $u_5v_5>1$ this implies 
$u_1v_1+ \cdots +u_4v_4\le u_5v_5-1$ which means $\langle \vect u,\vect v\rangle \le -1$. 

Furthormore, the proof that $\langle \vect u,\vect v\rangle \ne -1$ should be added:

Assume $\langle \vect u,\vect v\rangle=-1$. 
Since $\langle \vect u,\vect u\rangle=-1$ we have $\langle \vect u,\vect u-\vect v\rangle=0$. 
Lemma 9.1.1 (3) implies that $\vect u-\vect v$ is either space-like or equal to $\vect 0$. 
As $\langle \vect u-\vect v,\vect u-\vect v\rangle=0$, $\vect u-\vect v$ cannot be space-like. 
Therefore, $\vect u-\vect v=\vect 0$, which is a contradiction. 

\item Page 145. The 3rd line. $S^3_{infty}$ should be replaced by $S^3_{\infty}$. 

\item Page 145. The bijection in Theorem 9.3.2 can be considered as a modern version of the pentaspherical coordinates in [Dar]. 

\item Page 151 last two lines to Page 152 line. (1a). 
The book misses the description of the case when $n=1$. 

\smallskip
When $n=1$ $P^{\perp}\cap S^1_{\infty}=\emptyset$ since $P^{\perp}$ is a time-like line. 
In this case we may consider the ``base sphere" to be $\emptyset=\partial(P^{\perp}\cap\mathbb{H}^2)$, where $\partial\mathbb{H}^2=S^1_{\infty}$. 
Let us call $\mathcal{P}$ or the set of corresponding spheres of dimension $0$ a ``{\em space-like pencil}\,". 

\item Page 166. Three lines before {\bf Lemma 9.8.2}. 
$$\psi:\mathbb{R}^{n+1}_+\ni (\vect X,r)\vect{=}\varphi^{-1}\circ p^{-1}\left(S^{n-1}_r(\vect X)\right)\in\La\setminus ({\rm Span}\left<Q\right>)^{\perp}$$
should be replaced by
$$\psi:\mathbb{R}^{n+1}_+\ni (\vect X,r)\vect{\mapsto}\varphi^{-1}\circ p^{-1}\left(S^{n-1}_r(\vect X)\right)\in\La\setminus ({\rm Span}\left<Q\right>)^{\perp}.$$

\item Page 167. The 8th line. 

\begin{center}$\omegasub{\mathbb{R}^{n+1}_+}=\,${\footnotesize $\displaystyle \frac{\,1\,}{\vect{x_{n+2}}}$}$\,dx_1\w \cdots \w dx_{n+1}$ \end{center}

should be replaced by 

\begin{center}$\omegasub{\mathbb{R}^{n+1}_+}=\,${\footnotesize $\displaystyle \vect{-}\frac{\,1\,}{\,\vect{r^{n+1}}\,}$}$dX_1\w \cdots \w dX_n\vect{\w dr}$ \end{center}

\item Page 167. The 16th line. 

$$\widetilde g:\mathbb{R}^{n+1}_+\ni (\vect X, r)\mapsto \left(\frac{\vect X}{|\vect X|^2-r^2}\,, \,
\frac{r}{|\vect X|^2-r^2}\right)\in\mathbb{R}^{n+1}_+\,.$$

should be replaced by 

$$\widetilde g:\mathbb{R}^{n+1}_+\ni (\vect X, r)\mapsto \left(\frac{\vect X}{|\vect X|^2-r^2}\,, \,
\frac{r}{\,\vect{\vert}|\vect X|^2-r^2\vect{\vert}\,}\right)\in\mathbb{R}^{n+1}_+\,.$$

\item Page 175. {\bf Definition 10.1} (2) should be replaced by 

``(2) An oriented $2$-sphere $\varSigma$ 
is called a {\bf non-trivial sphere in the strict sense} for a knot $K$ 
if each connected component of $\mathbb{R}^3\setminus \varSigma$ ($S^3\setminus \varSigma$) 
contains at least 2 connected components of $K\setminus (K\cap\varSigma)$. "

\item Page 184. {\bf Definition 10.6}. (1) 

$$
{\mathcal C}c^{(4)}(K)=\left\{
(s,s,s,s)\in \Delta^{(4)}
\left| 
\begin{array}{l}
\mbox{\rm The osculating circle of $K$ at $f(s)$}\\
\mbox{\rm is tangent to $K$ to the {\bf fourth} order}
\end{array}
\right. \right\}.\index{${\mathcal C}c^{(4)}(K)$}
$$

should be replaced by 

$$
{\mathcal C}c^{(4)}(K)=\left\{
(s,s,s,s)\in \Delta^{(4)}
\left| 
\begin{array}{l}
\mbox{\rm The osculating circle of $K$ at $f(s)$}\\
\mbox{\rm is tangent to $K$ to the {\bf third} order}
\end{array}
\right. \right\}.\index{${\mathcal C}c^{(4)}(K)$}
$$

\item Page 184. {\bf Proposition 10.3.3}. (1) 

``The order of tangency of the osculating circle at $f(s)$ is 
exactly $\vect 3$ if and only if 
$f(s)\w f^{\p}(s)\w f^{\p\p}(s)\w f^{\p\p\p}(s)\ne \vect 0$. ''

should be replaced by 

``The order of tangency of the osculating circle at $f(s)$ is 
exactly $\vect 2$ if and only if 
$f(s)\w f^{\p}(s)\w f^{\p\p}(s)\w f^{\p\p\p}(s)\ne \vect 0$. ''

\item Page 185. The 6th line and the 10th line from the bottom (in the {\it Proof} of Proposition 10.3.3 (1)). 

``(1) Suppose the osculating circle $C$ is tangent to the knot at $f(s)$ 
to the {\bf fourth} order. $\cdots$ 

$\cdots$ and hence the osculating circle $C$ is tangent to the knot at $f(s)$ 
to the {\bf fourth} order. ''

should be replaced by 

``(1) Suppose the osculating circle $C$ is tangent to the knot at $f(s)$ 
to the {\bf third} order. $\cdots$ 

$\cdots$ and hence the osculating circle $C$ is tangent to the knot at $f(s)$ 
to the {\bf third} order. ''

\item Page 191. The 7th line. 

\[(I\times I)^{\ast}\la _{\Re}=\displaystyle{\la _{\Re}+\frac{\,1\,}2d\log (\vect{|u|^2})}\]

should be replaced by 

\[(I\times I)^{\ast}\la _{\Re}=\displaystyle{\la _{\Re}+\frac{\,1\,}2d\log (\vect{u^2+v^2}
)}\]

\item Page 193. The 4th line in the Proof of Lemma 11.2.2. 
`{\bf bas}' should be replaced by `{\bf basis}'. 

\item Page 193. The 5th line in Lemma 11.2.3. 
`{\bf this} local coordinate' should be replaced by `{\bf these} local coordinate{\bf s}'. 

\item Page 198. Theorem 11.2.7. As a corollary of this theorem, we have

\smallskip
{\bf Corollary} 
The pull-back $\omega ={\psi}^{\ast}\omegasub{0}$ of the canonical symplectic form $\omegasub{0}$ of $T^{\ast}S^n$ by $\psi :S^n\times S^n\setminus\Delta \to T^{\ast}S^n$ is invariant under any diagonal action of the M\"obius group ${\mathcal M}$ on $S^n\times S^n\setminus\Delta $. 

\item Page 200. The condition (ii) of Theorem 11.2.9 is not necessary. 

\item Page 203. A remark on Definition 11.4: 

The real part of the infinitesimal cross ratio of a knot $K$ is a smooth $2$-form on $K\times K\setminus\triangle$, but it is not the case with the imaginary part. 
Since the conformal angle $\theta_K(x,y)$ is not a smooth function of $x$ and $y$ (see Figure 10.1 on page 183), {\bf the imaginary part of the infinitesimal cross ratio may have singularity} at a pair of points $(x,y)\in K\times K \setminus\triangle$ where the conformal angle $\theta_K(x,y)$ vanishes. 

\item Page 206. The 2nd and 3rd lines
$$
\big((T_0\circ f)^{\ast}dx_3\big)(s_0, t_0)
=\left(\frac{d}{ds}(T_0\circ f)_3\right)(s_0), 
$$
$$
\big((T_0\circ f)^{\ast}dy_3\big)(s_0, t_0)
=\left(\frac{d}{dt}(T_0\circ f)_3\right)(t_0)  {}
$$
should be replaced by 
$$
\big((T_0\circ f)^{\ast}dx_3\big)(s_0, t_0)
=\left(\frac{d}{ds}(T_0\circ f)_3\right)(s_0)\vect{=0}, 
$$
$$
\big((T_0\circ f)^{\ast}dy_3\big)(s_0, t_0)
=\left(\frac{d}{dt}(T_0\circ f)_3\right)(t_0)\vect{=0}  {}
$$
i.e. $=0$ should be added at the end of the both formulae. 

\item Page 211. The 2nd and 3rd lines
$$\begin{array}{rcl}
\frac{\,1\,}2E^{(2)}_{\circ}(L_r)=E^{(2),mut}_{\circ}(L_r)
&=&\displaystyle{\int_0^{2\pi}\int_0^{2\pi}\frac{dsdt}{2\vect{+2}\cos (t-s)} }\\[3mm]
&=&\displaystyle{4\pi\int_0^{\frac{\,\pi\,}2}\frac{d\xi}{1-r\cos 2\xi} 
\hskip 1.0cm (\vect{\xi=\tan \eta}) }
\end{array}
$$
should be replaced by 
$$\begin{array}{rcl}
\frac{\,1\,}2E^{(2)}_{\circ}(L_r)=E^{(2),mut}_{\circ}(L_r)
&=&\displaystyle{\int_0^{2\pi}\int_0^{2\pi}\frac{dsdt}{2\vect{-2r}\cos (t-s)} }\\[3mm]
&=&\displaystyle{4\pi\int_0^{\frac{\,\pi\,}2}\frac{d\xi}{1-r\cos 2\xi} 
\hskip 1.0cm (\vect{\eta=\tan \xi}) }
\end{array}
$$

\item Page 218. The 3rd line in the Proof of Lemma 12.3.2 
$$
V_{\sin\theta}(K;f(\pm t)), V_{\sin\theta}(K;f(\de\pm t)) \vect{\ge}
\frac{\,1\,}{100}\cdot\frac{\,1\,}{d+t}. 
$$
should be replaced by 
$$
V_{\sin\theta}(K;f(\pm t)), V_{\sin\theta}(K;f(\de\pm t)) \vect{\le}
\frac{\,1\,}{100}\cdot\frac{\,1\,}{d+t}. 
$$
i.e. the inequality should be reversed. 

\item Page 218. The 8th line in the Proof of Lemma 12.3.2, i.e. just above Figure 12.2. 
$$
\lim_{K\setminus\{x_+\}\ni y\to x_+}\pi_{x_+}(\vect{y})=\hat x_+. 
$$
should be replaced by 
$$
\lim_{K\setminus\{x_+\}\ni y\to x_+}\pi_{x_+}(\vect{I_{x_+}(y)})=\hat x_+. 
$$

\item Page 218. The right picture of Figure 12.2 should be replaced by Figure \ref{Fig_12.2}. 
\begin{figure}[htbp]
\begin{center}
\includegraphics[width=.5\linewidth]{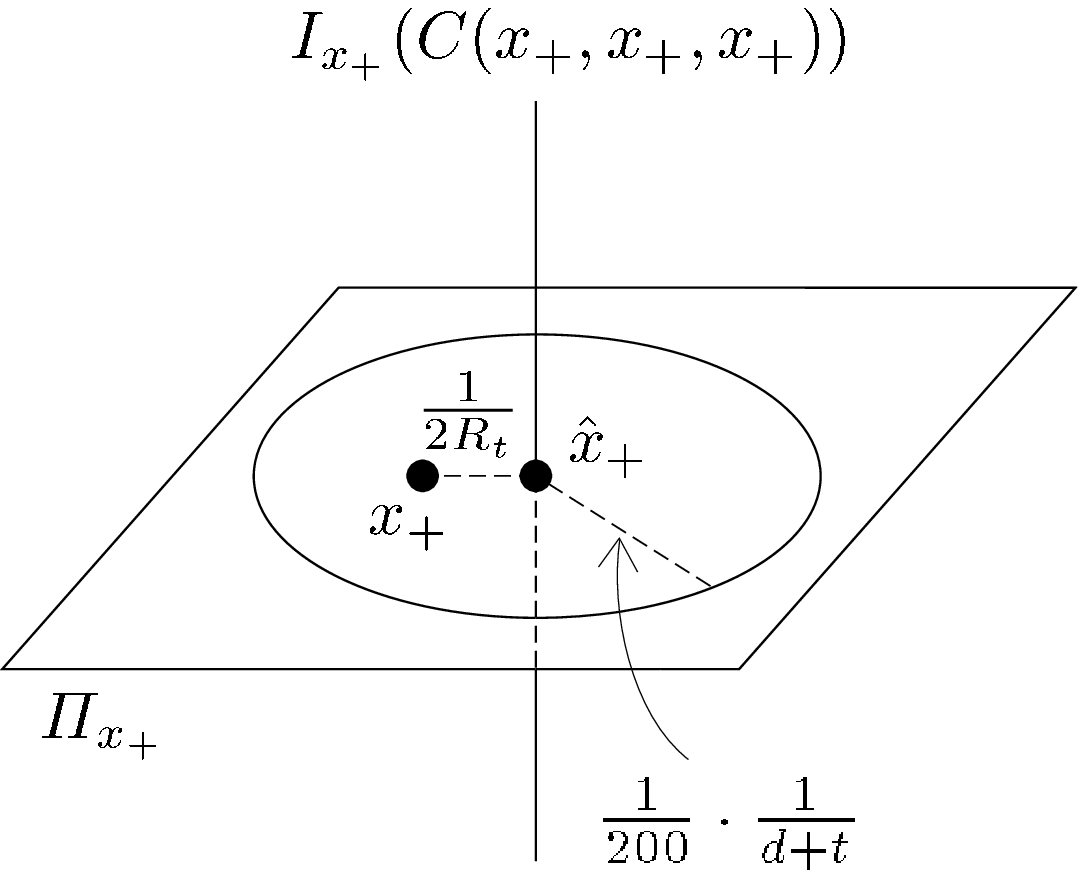}
\caption{The right picture of Figure 12.2.}
\label{Fig_12.2}
\end{center}
\end{figure}
The radius of the circle, which was $\frac1{400}\cdot\frac1{d+t}$ in the book, should be replaced by $\frac1{200}\cdot\frac1{d+t}$. 
Also, the point $x_+$ should be closer to the center of the circle, $\hat x_+$, as $|x_+-\hat x_+|\le \frac{\,1\,}{400}\cdot\frac{\,1\,}{\,d+t\,}$, which is a half of the radius. 

\item Page 218. The 3rd line from the bottom. 

$\pi_{x_+}(\widetilde{K}_{x_+})$ lies inside the circle on $\varPi_{x_+}$ 
with center $\vect{x_+}$ and radius $1/(200(d+t))$. 

should be replaced by 

$\pi_{x_+}(\widetilde{K}_{x_+})$ lies inside the circle on $\varPi_{x_+}$ 
with center $\vect{\hat x_+}$ and radius $1/(200(d+t))$. 

\item Page 219. The 8th line (just above Figure 12.3). 

$N_t$ whose meridian disc has radius $(\vect{400}(d+t))/3$, as illustrated in Figure 12.3. 

should be replaced by 

$N_t$ whose meridian disc has radius $(\vect{200}(d+t))/3$, as illustrated in Figure 12.3. 

\item Page 219. Figure 12.3 should be replaced by Figure \ref{Fig_12.3}. 
\begin{figure}[htbp]
\begin{center}
\includegraphics[width=.9\linewidth]{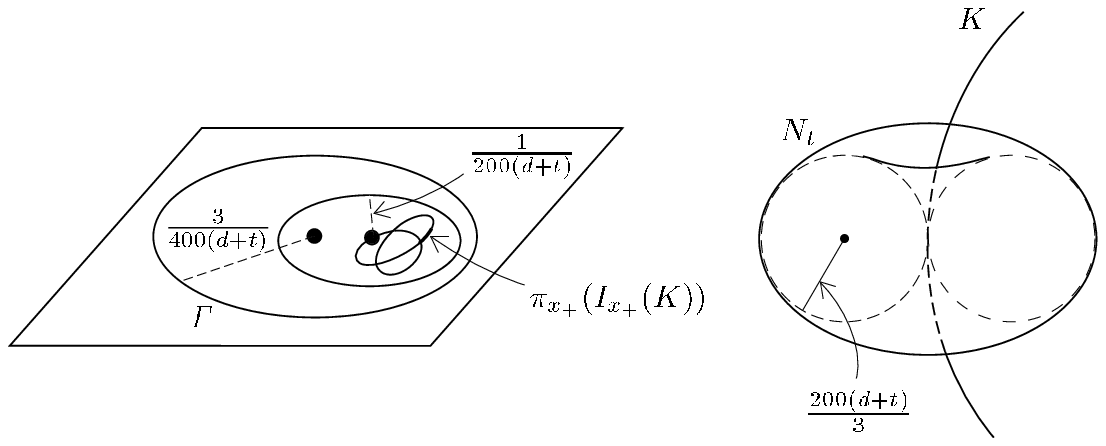}
\caption{Figure 12.3.}
\label{Fig_12.3}
\end{center}
\end{figure}
Namely, the radius of the inner circle of the left picture, which was $\displaystyle \frac1{\smallvect{400}(d+t)}$ in the book, should be $\displaystyle \frac1{\smallvect{200}(d+t)}$, and the radius of the `degenerate solid torus in the right picture, which was $\displaystyle \frac{\smallvect{400}(d+t)}3$ in the book, should be $\displaystyle \frac{\smallvect{200}(d+t)}3$. 

\item Page 219. The 7th line from the bottom 

$$|f(t)-f(-t)|\le 2t\ll\frac{\vect{400}}3(d+t)$$

should be replaced by 

$$|f(t)-f(-t)|\le 2t\ll\frac{\vect{200}}3(d+t)$$

\item Pae 220. The 2nd line from the bottom 

$$|f(0)-\vect{g(t)}|\ge\int_0^t(1-\kappa  s)ds=t\left(1-\frac{\,\kappa \,}2t\right)>d, $$

should be replaced by 

$$|f(0)-\vect{f(t)}|\ge\int_0^t(1-\kappa  s)ds=t\left(1-\frac{\,\kappa \,}2t\right)>d, $$

\item Page 226. The last three lines of Remark 13.2.2 should be replaced by 

This is because $1\le {}_{n}C_{2}\le{}_{2n}C_{4}$ when $n\ge 2$. 

\item Page 229. Line 6 up. 

``subarc of $C_3(r)$ between $Q_{12}(r)$ and $P_{{\bf 2}3}(r)$ "

should be replaced by 

``subarc of $C_3(r)$ between $Q_{12}(r)$ and $P_{{\bf 1}3}(r)$ "

\item Page 229. The last line. Page 230. The second line. 

The domain of integration ``$\displaystyle \left[\frac{\,1\,}{2\sin\theta_K}, \infty\right)$" should be replaced by ``$\displaystyle \left[\frac{\,1\,}{2\sin\theta}, \infty\right)$"

\item Page 231. Lines 3-5. 
``$\cdots +o(|z|)$" should be replaced by ``$\cdots +O(|z|)$". 

\item Page 236. The 13th line around the middle of the page 
$$
\int _{\varSigma _r(X,Y,Z)\cap K_d\ne\smallvect{\varphi}}  \sharp(\varSigma _r(X,Y,Z)\cap K_d)\,dXdYdZ
=2\pi r^2 \left(\bar K_d\right), 
$$
should be replaced by 
$$
\int _{\varSigma _r(X,Y,Z)\cap K_d\ne\smallvect{\emptyset}}  \sharp(\varSigma _r(X,Y,Z)\cap K_d)\,dXdYdZ
=2\pi r^2 \left(\bar K_d\right), 
$$

\item Page 236. The 5th line from the bottom
$$E_{{\rm mnts}}(K_d)=\int _{\mathcal{S}(K_d)}{}
C^{\sharp(\varSigma _r(X,Y,Z)\cap K_d)}_2
\cdot\frac{\,1\,}{\,r^4\,}\,dXdYdZdr 
$$
should be replaced by 
$$
E_{{\rm mnts}}(K_d)=\int _{\mathcal{S}(K_d)}{}
{}_{\smallvect{\frac12}\sharp(\varSigma _r(X,Y,Z)\cap K_d)}C_2
\cdot\frac{\,1\,}{\,r^4\,}\,dXdYdZdr 
$$

This implies the 4th line from the bottom because ${}_mC_2\ge m-1$ for $m\ge 0$. 

\item Page 256. The 8th line (just before Example A.1) should be replaced by 

We show that $\vect{I_{tv}}$ can detect the unknot. 

Namely, it is a conjecture whether $|csl|$ can detect the unknot or not. 

\item Page 273. 

[CKS2] J. Cantarella, R.~B. Kusner, and
 J.~M. Sullivan, 
{\em On the minimum ropelength of knots and links,} Preprint. 

should be updated to 

[CKS2] J. Cantarella, R.~B. Kusner, and
 J.~M. Sullivan, 
{\em On the minimum ropelength of knots and links,} Invent. Math. {\bf 150}  (2002),  no. 2, 257--286. 

[CKS] should be removed as it is same as [CKS1].

\item Page 276. 

[GLP] M. Gromov, J. Lafontaine, and P. {\bf Pnau}, 
{\em Structures m\'etriques pour les vari\'et\'es riemanniennes,} 
Cedic/Fernand Nathan, Paris, 1981.

should be replaced by 

[GLP]{Gr-La-Pn} M. Gromov, J. Lafontaine, and P. {\bf Pansu}, 
{\em Structures m\'etriques pour les vari\'et\'es riemanniennes,} 
Cedic/Fernand Nathan, Paris, 1981.

\item Page 279. 

[Lin] X.-S. Lin, 
{\em Knot energies and knot invariants,} {\bf J. Differential Geom. {\bf 44} (1996), 74\,--\,95}.

should be replaced by 

[Lin] X.-S. Lin, 
{\em Knot energies and knot invariants. Knot theory and its applications,} 
{\bf Chaos Solitons Fractals {\bf 9} (1998) 645\,--\,655}

\item Page 285, Index. The ``average crossing number" also appears 
on page 44. 


  

\end{itemize}

\bigskip
\noindent
{\large{\bf Acknowledgement}}. The author thanks Rob Kusner for pointing out a couple of mistakes. He also thanks Paul Feehan for helpful comments. 







\end{document}